 \documentclass[final,5p,times,twocolumn]{elsarticle}

\usepackage{amssymb}
\usepackage{txfonts}

\usepackage{graphicx}
\usepackage{subcaption}
\usepackage{float}

\usepackage{amsmath}
\usepackage{CJKutf8}
\usepackage{tabularx}
\usepackage{hyperref}

\usepackage{xcolor}
\hypersetup{
	colorlinks=true,
	urlcolor=blue, 
	citecolor=blue 
}
\usepackage{natbib} 

\usepackage{booktabs}

\journal{Arxiv}

\begin{document}
	
	\newtheorem{theorem}{Theorem}[section] 

	\newtheorem{lemma}{Lemma}[section]
	
	\begin{frontmatter}
		
		
		
		\title{3-path-connectivity of Cayley graphs generated by wheel graphs} 
		\author[1]{Yi-Lu Luo}
		\author[1]{Yun-Ping Deng \corref{cor1}}
		\ead{dyp612@hotmail.com}
		\author[1]{Yuan Sun}
		\cortext[cor1]{Corresponding author}
	
		\affiliation[1]{organization={Department of Mathematics, Shanghai University of Electric Power}, 
			addressline={2103 Pingliang Road},
			postcode={200090}, 
			city={Shanghai}, 
			country={PR China}}
		
		\begin{abstract}
			Let $G = (V(G), E(G))$ be a simple connected graph and $\Omega$ a subset of $ V(G)$ with $|\Omega|\geq2$. An  $\Omega$-path in $G$ is a path that connects all vertices of $\Omega$. Two $\Omega$-paths $P_i$ and $P_j$ are said to be internally disjoint if $V(P_i)\cap V(P_j)=\Omega$ and $E(P_i)\cap E(P_j)=\emptyset$. Denote $\pi_G(\Omega)$ by the maximum number of internally disjoint $\Omega$-paths in $G$. For an integer $k\geq2$, the $k$-path-connectivity $\pi_k(G)$ of $G$ is defined as  $\min\{\pi_G(\Omega)\mid\Omega\subseteq V(G)$ and $|\Omega|=k\}$. Let $CW_n$ denote the Cayley graph generated by the $n$-vertex wheel graph. In this paper, we investigate the $3$-path-connectivity of $CW_n$ and prove that  $\pi_3(CW_n)=\lfloor\frac{6n-9}4\rfloor$ for all $n\geq4$.
		\end{abstract}
		
		
		
		
	\end{frontmatter}
	
		
		\section{Introduction}\label{sec1}

	An interconnection network is usually modeled as a graph, and connectivity serves as a crucial parameter for assessing the reliability and fault tolerance of the graph. 
	Let $G$ be a simple connected graph with vertex set $V(G)$ and edge set $E(G)$.
	For any $2$-subset $\{u, v\} \subseteq V(G)$, let $\kappa_G(u, v)$ denote the maximum number of
	internally disjoint paths connecting vertices $u$ and $v$ in $G$.
	The vertex connectivity $\kappa(G)$ of $G$ is defined as $\min\{\kappa_G(u, v)\mid\{u, v\} \subseteq V(G)\}$.
	As an extension of vertex connectivity, Hager \cite{Hager86} proposed the concept of $k$-path-connectivity, which generalizes the connection of two vertices to that of $k$ vertices (where $k\geq2$).
	Let $\Omega\subseteq V(G)$ with $|\Omega|\geq2$. A path that connects all vertices of $\Omega$ is called an $\Omega$-path of $G$. Two $\Omega$-paths $P_i$ and $P_j$ are said to be internally disjoint if $V(P_i)\cap V(P_j)=\Omega$ and $E(P_i)\cap E(P_j)=\emptyset$. 
	Let $\pi_G(\Omega)$ denote the maximum number of internally disjoint $\Omega$-paths in $G$. 
	For an integer $k\geq2$, the $k$-path-connectivity $\pi_k(G)$ of $G$ is defined as $\min\{\pi_G(\Omega)\mid\Omega\subseteq V(G)$ and $|\Omega|=k\}$.
	The path-connectivity has attracted much attention since it was proposed. Li et al. \cite{Li21} proved that for any $k \geq1$, deciding whether $\pi_G(T)\geq k$ with $T \subseteq V(G)$ is NP-complete. Zhu et al. \cite{Zhu23} established an upper bound for the $3$-path connectivity of regular connected graphs, and determined the $3$-path connectivity of the $k$-ary $n$-cube. Jin et al. \cite{Jin23} characterized the $3$-path connectivity of Cayley graphs generated by transposition trees.
	Recently, Wang et al. \cite{wang25} determined the $3$-path connectivity of pancake graphs. 
	For more results about the $k$ path-connectivity of graphs, we refer the reader to \cite{Li24,li2025k,luo20253,Mao16,Zhu22,zhu2024reliability}. 
	

	Cayley graphs possess some properties that are highly desirable in the context of interconnection networks \cite{heydemann1997cayley,Lakshmivarahan93}.
	Let $X$ be a finite group and $S$ an inverse-closed subset of $X \backslash \{e\}$, where $e$ is the identity element of $X$.
	The Cayley graph $Cay(X, S)$  is the graph with vertex set $X$ and edge set $\{(g, gs) \,|\, g\in X, s\in S\}$. 
	Let $S_n$ be the symmetric group of degree $n$. 
	Define the generating sets   $S_1=\{(1~j)\mid2\leqslant j\leqslant n\}\cup\{(j~j{+}1)\mid2\leqslant j\leqslant n-1\}$ and $S_2 = \{(1~j)\mid2\leqslant j\leqslant n\}\cup\{(j~j{+}1)\mid2\leqslant j\leqslant n-1\}\cup\{(2~n)\}$. The bubble-sort star graph $BS_n$ \cite{Chou96} and the Cayley graphs $CW_n$ generated by
	wheel graphs \cite{luo2020kind} are then  defined as $Cay(S_n, S_1)$ and $Cay(S_n, S_2)$, respectively. 
	$CW_n$ represents a typical topological structure in interconnection networks, and its connectivity properties have been extensively studied  \cite{feng2019nature,feng20202,feng20,feng20212,feng2022diagnosability,tu2017kind,zhao2019generalized}. In this paper, we focus on the $3$-path-connectivity of $CW_n$ and prove that $\pi_3(CW_n)=\lfloor\frac{6n-9}4\rfloor$ for all $n\geq4$.

	\section{Preliminaries}
	
	Let $G$ be a simple connected graph with vertex set $V(G)$ and edge set $E(G)$.
	Table \ref{table1} presents a comprehensive list of notations used throughout the paper. For terminology and notations not defined here, we refer the reader to \cite{Bondy08}.
	For any $i \in [n]$, let $CW_n^{i}$ be the subgraph of $CW_n$ induced by the set $\{ \sigma \in S_n \mid \sigma(n)=i\}$. 
	Obviously, 
	$CW_n$ can be decomposed into $n$ disjoint 
	subgraphs $CW_n^i \,(i\in[n])$, each isomorphic to $BS_{n-1}$, i.e. a copy of $BS_{n-1}$.
	 We denote this decomposition by ${CW_n}=CW_{n}^{{1}}\oplus CW_{n}^{{2}}\oplus\cdots\oplus CW_{n}^{{n}}$.
	For any subset $\{i_1,i_2,\ldots,i_t\}\subseteq[n]$, let  $CW_{n}^{i_{1}}\oplus CW_{n}^{i_{2}}\oplus\cdots\oplus CW_{n}^{i_{t}}$ denote the subgraph of $CW_n$ induced by $\bigcup_{k=1}^{t} V(CW_n^{i_k})$.

Next we shall state some facts on $BS_{n}$ and $CW_n$ without proof in the following lemma. Some of these facts may be found in \cite{Cai15,feng2019nature,feng20202}, and the others follow immediately from the definition of $CW_n$.
	
	\begin{lemma}\label{lem2.1} (see \cite{Cai15,feng2019nature,feng20202})
		For $n\geq 4$, $i,j\in [n]$ and $i\neq j$, the following properties hold:		
		
		(i) $BS_{n}$ is $(2n-3)$-regular, vertex transitive, and $\kappa(BS_{n})=2n-3$ for $n\geq2$.
		

		(ii) For any $2$-subset $\{u, v\}\subseteq V(CW_n)$, $\left|CN_{CW_n}(u,v)\right|\leq 3$.
		
		(iii) $\left|E_{i, j}(CW_n)\right|=3(n-2)!$, where $E_{i, j}(CW_n)$ is the set of edges between $CW_n^{i}$ and $CW_n^{j}$.

		(iv) For any $v\in V(CW_n^{i})$, $v$ has exactly $2n-5$ neighbors in $CW_n^{i}$ and three neighbors outside $CW_n^{i}$. These three outside neighbors belong to three distinct copies and are denoted as:
		 $v^+ = v(1~n)$, $v^- = v(n{-}1~n)$, $v^* = v(2~n)$.
		Clearly, $\{u^+,u^-,u^*\}\cap\{v^+,v^-,v^*\}=\emptyset$ for any distinct vertices $u, v\in V(CW_n^{i})$.

	\end{lemma}
	
	\begin{table}[htbp]
		\centering
		\caption{Notation and its meaning}
		\renewcommand{\arraystretch}{1.4} 
		\begin{tabularx}{\columnwidth}{@{}lX@{}}
			\toprule
			Notation 
			& Meaning \\ \midrule
			$\kappa_G(u, v)$ & The maximum number of internally disjoint paths connecting vertices $u$ and $v$ in $G$ \\
			$G\setminus H$ & Subgraph of $G$ induced by $V(G)-V(H)$ \\
			$N_G(v)$ & The neighbors of the vertex $v$ in $G$ \\
			$CN_G(v_1, \ldots, v_k)$ & $N_G(v_1)\cap\cdots\cap N_G(v_k)$ \\
			$d_G(v)$ & The degree of the vertex $v$ in $G$ \\
			$[n]$ & $\{1,2,\ldots, n\}$ \\
			$[m,n]$ & $\{m,m+1,\ldots, n\}$ \\
			$d(P)$ & The length of path $P$ \\
			$d_P(u,v)$ & The distance between two vertices $u$ and $v$ in path $P$ \\
			$P(u,v)$ & The subpath from $u$ to $v$ in path $P$ \\
			$W.l.o.g.$ & Without loss of generality \\
			\bottomrule
		\end{tabularx}
		\label{table1}
	\end{table}
	
	\begin{lemma}\label{lem2.2} (see \cite{Bondy08})
		Let $G$ be a $k$-connected graph, let $x$ be a vertex of $G$, and let $Y \subseteq V(G) \setminus \{x\}$ be a set of at least $k$ vertices of $G$. Then there exists a $k$-fan in $G$ from $x$ to $Y$.
	\end{lemma}
	
	\begin{lemma}\label{lem2.3} (see \cite{Bondy08})
		Let $G$ be a $k$-connected graph, and let $X$ and $Y$ be subsets of $V(G)$ of cardinality at least $k$. Then there exists in $G$ a family of $k$ pairwise disjoint $(X,Y)$-paths.
	\end{lemma}

	\begin{lemma}\label{lem2.4} (see \cite{luo20253})
		For any $n\geq 3$, it holds that $max\left \{ \left | CN_{BS_n}( u,v,w) \right | \mid \{ u,v,w\} \subseteq V(BS_n) \right \} =3$.
	\end{lemma}

	\begin{lemma}\label{lem2.5} 
			For any $n\geq 4$, it holds that  $max\left \{\, \left | CN_{CW_n}( u,v,w) \right | \mid \{ u,v,w\} \subseteq V(CW_n) \,\right  \} =3$.
		\\
		\\
		\noindent\textbf{Proof.} 
	Since $BS_3$ is a subgraph of $CW_n$ for any $n\geq4$, it follows from Lemma \ref{lem2.4} that  $max\left \{ \left | CN_{CW_n}( u,v,w) \right | \mid \{ u,v,w\} \subseteq V(CW_n)\right  \} \geq3$.
		
		On the other hand, by Lemma \ref{lem2.1} (ii), $\left|CN_{CW_n}(u,v)\right|\leq3$, and thus $\left|CN_{CW_n}(u,v,w)\right|$ $\leq3$ for any $u,v,w\in V(CW_n)$, that is, $max\left \{\left|CN_{CW_n}(u,v,w)\right|\mid\{u,v,w\}\subseteq V(CW_n)\right \} {\leq}3$.

		Thus, $max\left \{ \left | CN_{CW_n}( u,v,w) \right | \mid \{ u,v,w\} \subseteq V(CW_n) \right \}=3$. $\square$
	\end{lemma}


	\begin{lemma}\label{lem2.6}
	Let $\widehat{CW_n}=CW_{n}^{i_{1}}\oplus CW_{n}^{i_{2}}\oplus\cdots\oplus CW_{n}^{i_{t}}$ for any $\{i_1,\ldots,i_t\}\subseteq[n]$ with $n\geq4$.
		Then $\kappa(\widehat{CW_n})\geq2n-5$.
		\\
		\\
		\noindent\textbf{Proof.}
			When $t=1$, since $CW_{n}^{i_1}$ is isomorphic to $BS_ {n-1}$, it follows from Lemma \ref{lem2.1} (i) that
	 $\kappa(\widehat{CW_{n}})=\kappa(CW_{n}^{i_1})=\kappa(BS_{n-1})=2(n-1)-3=2n-5$.
		
		When $t \geq2$, let $v_1$ and $v_2$ be any two distinct vertices in $V(\widehat {CW_n})$, next we prove that $\kappa_{\widehat {CW_n}}(v_1,v_2)\geq 2n-5$ by considering the following two cases, from which we can conclude that $\kappa(\widehat{CW_n})\geq2n-5$.
		
		\textbf{Case 1.} $v_1 ,v_2$ belong to the same copy $CW_{n}^i$.
		
		Since $CW_{n}^i$ is isomorphic to $BS_ {n-1}$, it follows from Lemma \ref{lem2.1} (i) that $\kappa(CW_{n}^i)=2n - 5$, and thus  $\kappa_{\widehat {CW_n}}(v_1,v_2)\geq 2n-5$.
		
		\begin{figure}[h]
			\centering
			\includegraphics[width=0.9\columnwidth]{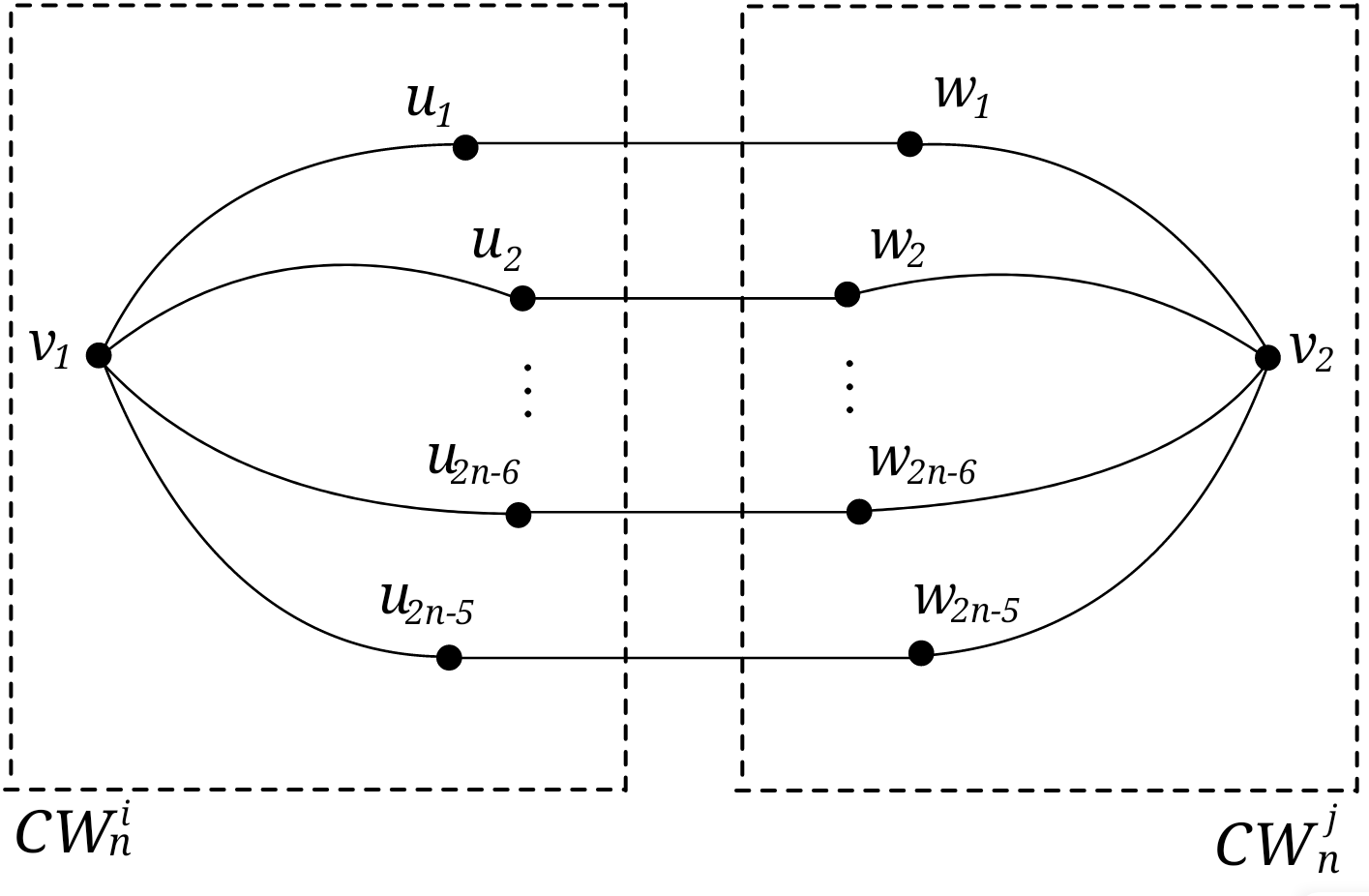} 
			\caption{Illustrations of Lemma 2.6 Case 2}
			\label{fig1}
		\end{figure}
		
		\textbf{Case 2.} $v_1$ and $v_2$ belong to two different copies $CW_{n}^i$ and $CW_{n}^j$, respectively.
		
		By Lemma \ref{lem2.1} (iii), $\left|E_{i,j}(CW_n)\right|=3(n-2)!$. Since $2n-5\leq3(n-2)!$ for any $n\geq4$,
	 referring to Figure \ref{fig1}, we can select $2n-5$ edges $(u_k,w_k)$ in $E_{i,j}(CW_n)$, where $u_k\in V(CW_{n}^i)\setminus\{v_1\}$ and $w_k\in V(CW_{n}^j)\setminus\{v_2\}$ with $1\leq k\leq 2n-5$.
		
		Let $X=\{u_1,u_2,\ldots,u_{2n-5}\}$ and $Y=\{w_1,w_2,\ldots,w_{2n-5}\}$.
		Since $\kappa(CW_{n}^i)=\kappa(CW_{n}^j)=2n-5$, it follows from Lemma \ref{lem2.2} that there are $2n-5$ internally disjoint $(v_1,X)$-paths $P_1, P_2,\ldots, P_{2n-5}$ in $CW_{n}^i$ and
		$2n-5$ internally disjoint $(v_2, Y)$-paths $R_1, R_2,\ldots, R_{2n-5}$ in $CW_{n}^j$. 
		Then we can obtain $2n-5$ internally disjoint $(v_1,v_2)$-paths: $v_1 P_i u_i w_i R_i v_2$, $i=1,2,\ldots 2n{-}5$.
		Hence, $\kappa_{\widehat {CW_n}}(v_1,v_2)\geq 2n-5$.
		 
	\end{lemma}

	\begin{lemma}\label{lem2.7}
		
		For $n\geq4$, $\kappa(CW_n\setminus{CW_n^k})\geq 2n-4$ with $k\in [n]$.
		\\
		\\
		\noindent\textbf{Proof.}
			Let $v_1 $ and $v_2$ be any two distinct vertices in $V(CW_n\setminus{CW_n^k})$. Next, we prove that $\kappa_{CW_n\setminus{CW_n^k}}(v_1,v_2)\geq 2n-4$ by considering the following two cases, which implies that $\kappa(CW_n\setminus{CW_n^k})\geq 2n-4$.
		
		\textbf{Case 1.} $v_1 ,v_2$ belong to the same copy $CW_{n}^i$.
		
		Since $CW_n^i$ is isomorphic to $BS_ {n-1}$, it follows from Lemma \ref{lem2.1} (i) that $\kappa(CW_{n}^i)=2n-5$, and thus there are $2n-5$ internally disjoint $(v_1, v_2)$-paths in $CW_{n}^i$. In addition, by Lemma \ref{lem2.1} (iv), $v_1$ and $v_2$ have at least one outside neighbor $v_1^{\prime},v_2^{\prime}\in V({CW_n}\setminus (CW_{n}^i\oplus CW_n^k))$, respectively.
		Since  $\kappa( {CW_n}\setminus (CW_{n}^i\oplus CW_n^k))\geq 1$, there exists a $(v_1^\prime,v_2^\prime)$-path in ${CW_n}\setminus (CW_{n}^i\oplus CW_n^k)$, which can derive another  $(v_1,v_2)$-path.
		Thus, $\kappa_{CW_n\setminus{CW_n^k}}(v_1,v_2)\geq 2n-4$.
		
		\textbf{Case 2.} $v_1$ and $v_2$ belong to two different copies $CW_{n}^i$ and $CW_{n}^j$, respectively.
		
			Similarly to the proof of Case 2 in Lemma \ref{lem2.6} , there are $2n-5$ internally disjoint $(v_1, v_2)$-paths in $CW_{n}^i\oplus CW_n^j$. In addition, by Lemma \ref{lem2.1} (iv), $v_1$ and $v_2$ have at least one outside neighbor $v_1^{\prime},v_2^{\prime}\in  V({CW_n}\setminus(CW_{n}^i\oplus CW_{n}^j\oplus CW_n^k))$, respectively. Since  $\kappa( {CW_n}\setminus (CW_{n}^i\oplus CW_n^j\oplus CW_n^k))\geq 1$, 
			there exists a $(v_1^\prime,v_2^\prime)$-path in ${CW_n}\setminus(CW_{n}^i\oplus CW_{n}^j\oplus CW_n^k)$, which can derive another  $(v_1,v_2)$-path.
		Thus, $\kappa_{CW_n\setminus{CW_n^k}}(v_1,v_2)\geq 2n-4$. $\square$
	\end{lemma}

	\begin{lemma}\label{lem2.8} (see \cite{luo20253})
		Let $\Omega=\{a,b,c\}$ be an arbitrary subset of $V(BS_n)$ with  $n=2d \, (where~d\geq2)$. Then in $BS_n$, there exist $2d-2$ $(a,b)$-paths, $2d-2$ $(a,c)$-paths and $2d-2$ $(b,c)$-paths. These $6d-6$ paths are internally disjoint, with none of their internal vertices belonging to $\Omega$.

	\end{lemma}

	\begin{lemma}\label{lem2.9} (see \cite{Zhu23})
		For any $k$-regular connected graph $G$, it holds that $\pi_3\left(G\right)\leq\lfloor\frac{3k-r}4\rfloor$,
		where $r=max \left \{ \left | CN_{G}( u,v,w) \right| \mid\{u,v,w\}\subseteq V({G})\right\}.$
	\end{lemma}

	\section{A structure connecting any three vertices in $V(CW_n)$}\label{sec3}
	In this section, considering the parity of $n$, we show that $CW_n$ has a structure connecting any three vertices in $V(CW_n)$.
	
	\begin{theorem}\label{the3.1} 
		For any $3$-subset $\Omega=\{a,b,c\}$ of $V(CW_n)$ with $n\geq4$, $CW_n$ contains the structure illustrated in Figure  \ref{fig2}.
		
		(1) If $n=2d\,(d\geq2)$, then there exist $2d-2$ $(a,b)$-paths, $2d-2$ $(a,c)$-paths and $2d-2$ $(b,c)$-paths. These $6d-6$ paths are internally disjoint, with none of their internal vertices belonging to $\Omega$.
		
		(2) If $n=2d+1\,(d\geq2)$, then there exist $2d-2$ $(a,b)$-paths, $2d$ $(a,c)$-paths and $2d$ $(b,c)$-paths. These $6d-2$ paths are internally disjoint, with none of their internal vertices belonging to $\Omega$.
		
	\begin{figure}[htbp]
		\centering
		\subfloat[$n=2d\,(d\geq2)$]{
			\includegraphics[width=0.48\columnwidth]{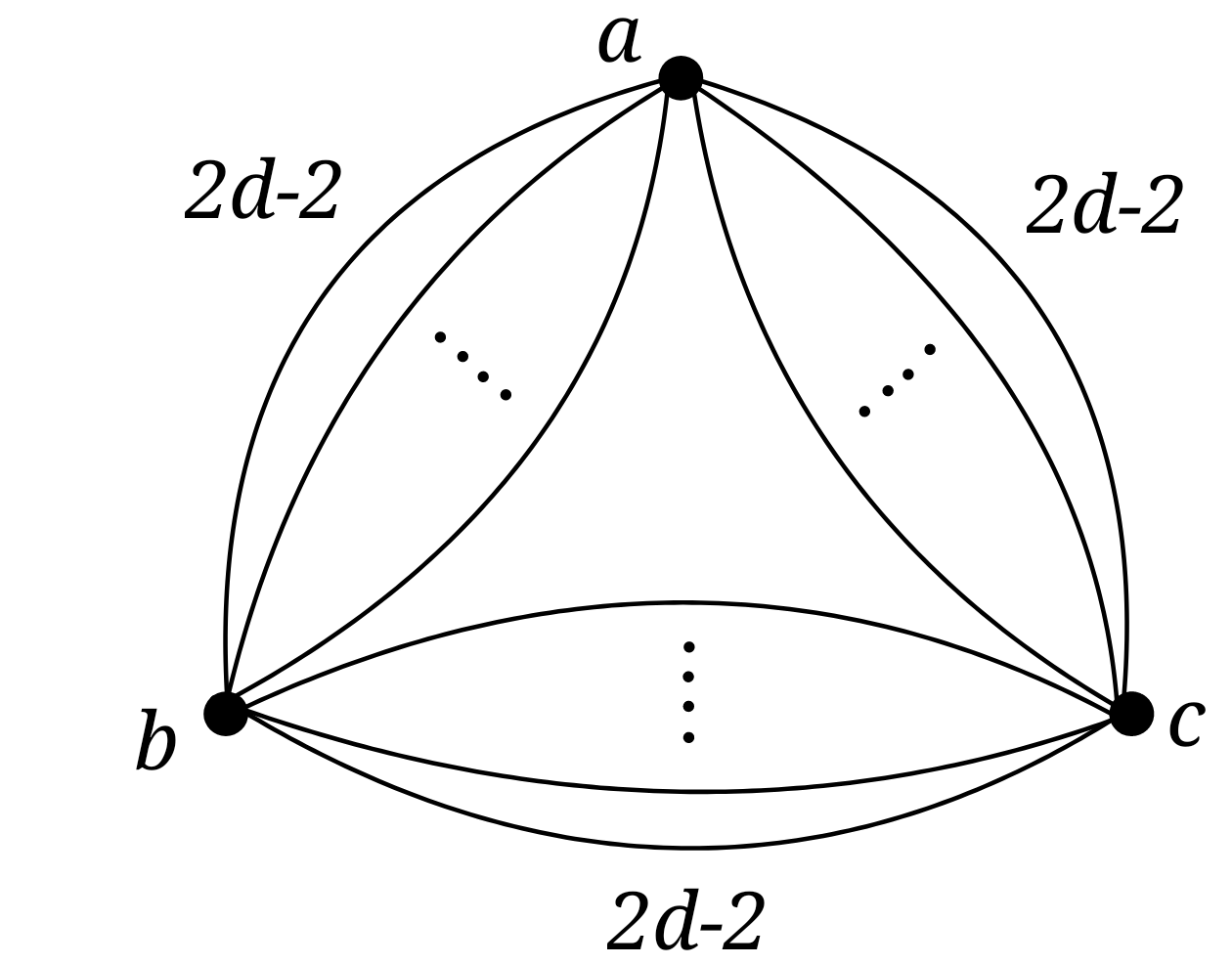}
		
		}
		\hfill
		\subfloat[$n=2d+1\,(d\geq2)$]{
			\includegraphics[width=0.48\columnwidth]{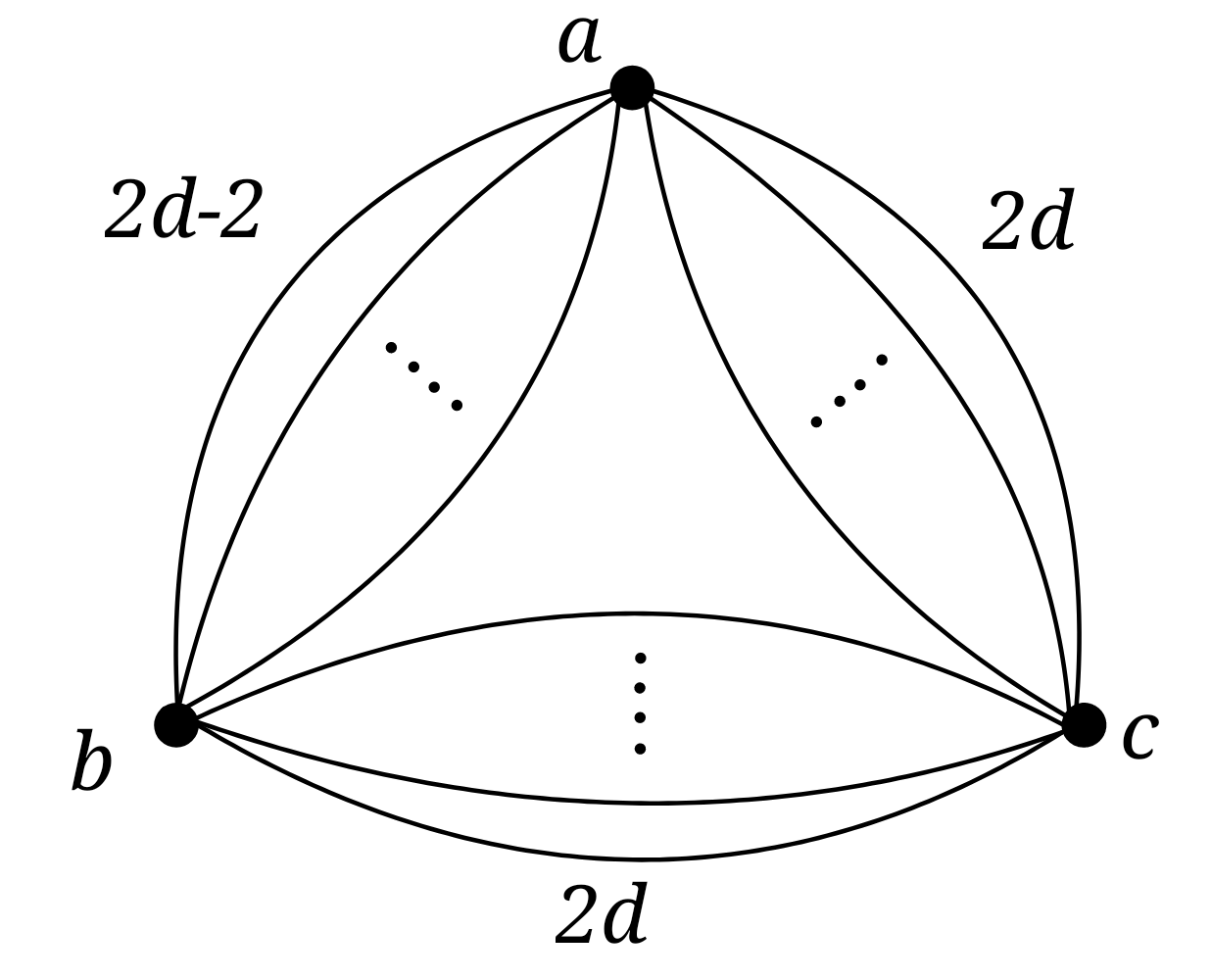}
	
		}
		\caption{Illustrations of Theorem 3.1}
		\label{fig2}
	\end{figure}

		\noindent\textbf{Proof.}
		For $n=2d \, (d\geq2)$, since $BS_n$ is a subgraph of $CW_n$, it follows from Lemma \ref{lem2.8} that the desired structure is obtained. Next we prove the situation when $n=2d+1 \, (d\geq2)$ by considering the following three cases:
		
		\textbf{Case 1.} $a,b$, and $c$ belong to the same copy.
		
		W.l.o.g., assume that $a,b,c\in CW_n^n$. 
		Since $CW_n^n$ is isomorphic to $BS_{n-1}$, it follows from Lemma \ref{lem2.8} that in $CW_n^n$, there are  $2d-2$ $(a,b)$-paths, $2d-2$ $(a,c)$-paths and $2d-2$ $(b,c)$-paths, which are denoted by $P_i,Q_i$ and $R_i$ $(i\in[2d-2])$, respectively.
		 Since $d_{CW_n^n}(a)=2n-5=4d-3$, there must exist a vertex $a^\prime\in N_{CW_n^n}(a)$ that does not appear in any of the above-mentioned $(a, b)$-paths $P_i$ or $(a, c)$-paths $Q_i$. Similarly, there must exist a vertex $b^\prime\in N_{CW_n^n}(b)$ not appearing in any $P_i$ or $R_i$, and a vertex $c^\prime\in N_{CW_n^n}(c)$  not appearing in any $Q_i$ or $R_i$ (for $i\in[2d-2]$). Next, we consider the following two subcases:
		
		\textbf{Subcase 1.1.} There exists a vertex in $\Omega$ that has an inside neighbor not appearing in any of the above-mentioned $6d-6$ paths.
		
		W.l.o.g., assume that $c$ has an inside neighbor $c^\prime$ not appearing in any of the above-mentioned $6d-6$ paths. Let $X=\{c^+,c^-,c^*,(c^\prime)^+\}$ and $Y=\{a^+,a^-,b^+,b^-\}$. 
		By  Lemma \ref{lem2.7}, $\kappa(CW_n\setminus CW_n^n)\geq2n-4\geq4$ for $n\geq 4$.  Consequently, Lemma \ref{lem2.3} implies the existence of four pairwise disjoint $(X,Y)$-paths in $CW_n\setminus CW_n^n$,  which in turn yield two additional $(a,c)$-paths and two additional $(b,c)$-paths. Combining these paths with $P_i, Q_i, R_i \, (i\in[2d-2])$,  the structure illustrated in Figure \ref{fig2} (b) is obtained.

		\textbf{Subcase 1.2.} All inside neighbors of each vertex in $\Omega$ belong to $P_i, Q_i$, or $R_i (i \in [2d-2])$.
		
		W.l.o.g., referring to Figure \ref{fig3} or \ref{fig4} , assume that $a^\prime$ belongs to a $(b,c)$-path $R_1$, $b^\prime$ belongs to an $(a,c)$-path $Q_1$, and $c^\prime$ belongs to an $(a,b)$-path $P_1$.
		Let $L = \{
		P_1(a,c'),\ P_1(c',b),\ Q_1(a,b'),\\ Q_1(b',c),\ R_1(b,a'), \ R_1(a',c)
		\}$.
		
		\begin{figure}[h]
			\centering
			\includegraphics[width=0.85\columnwidth]{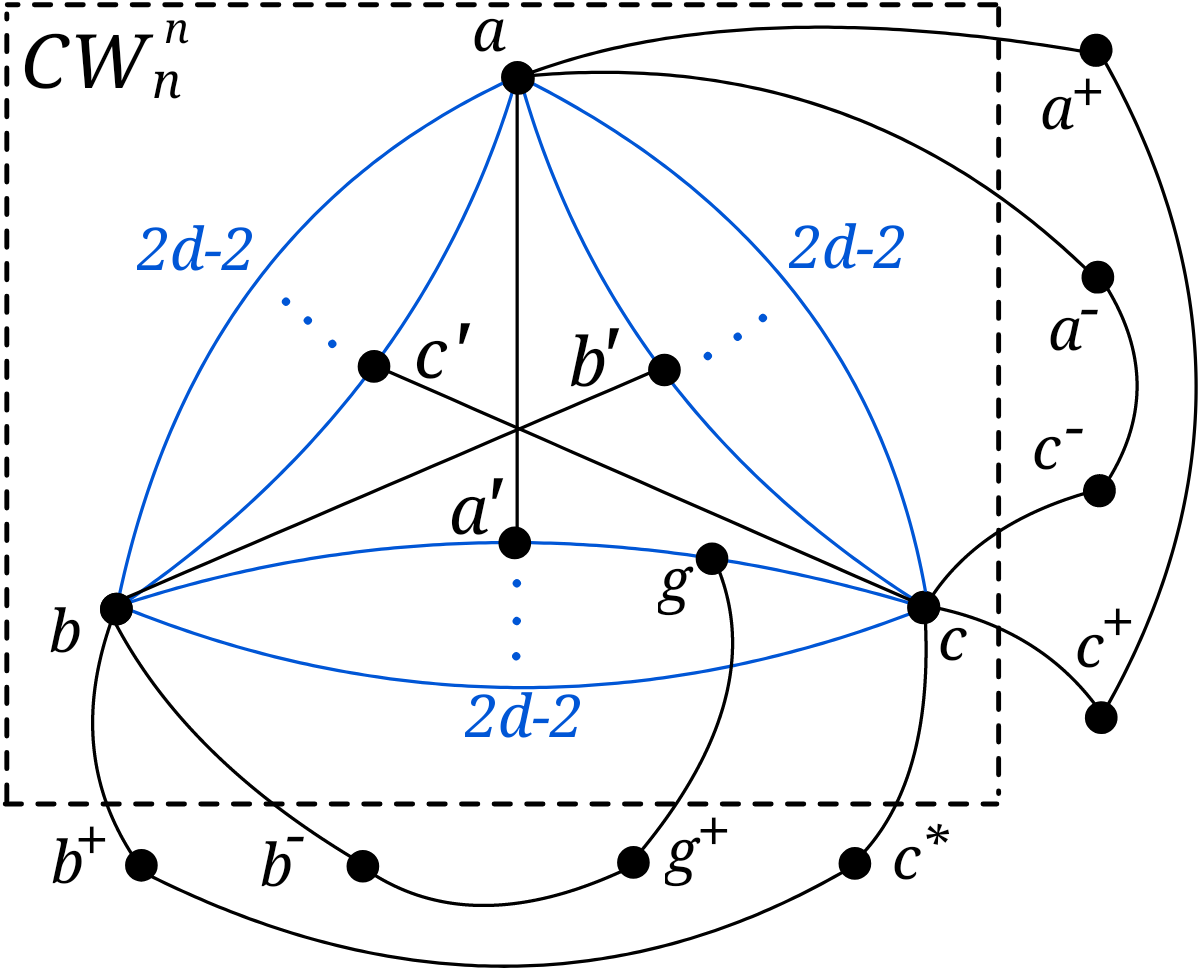} 
			\caption{Illustration of Subcase 1.2.1}
			\label{fig3}
		\end{figure}
		
		\textbf{Subcase 1.2.1.} 
		 There exists a path $P\in L$ such that $d(P)\geq2$.
		
		W.l.o.g., referring to Figure \ref{fig3}, assume that $d(R_1(a^\prime,c))=2$ and $N_{R_1}(c)=\{g\}$. 
		Let $X=\{c^+,c^-,c^*,g^+\}$ and $Y=\{a^+,a^-,b^+,b^-\}$.
		Similar to the proof of Subcase 1.1, there are four pairwise disjoint $(X, Y)$-paths in $CW_n\setminus CW_n^n$,  which in turn yield two additional $(a, c)$-paths denoted by $Q_{2d-1}, Q_{2d}$ and  two additional $(b, c)$-paths denoted by $R_{2d-1}, R_{2d}$. Let $P_1^* = aa^\prime \cup R_1(a^\prime,b) , Q_1^* = P_1(a,c^\prime) \cup c^\prime c$,  and
		$R_1^* = bb^\prime \cup Q_1(b^\prime,c)$ . Then, we have obtained the structure illustrated in Figure \ref{fig2} (b), which includes $2d-2$ $(a, b)$-paths: $P_1^*, P_2, \ldots , P_{2d-2}$,  $2d$ $(a, c)$-paths: $Q_1^*, Q_2,\ldots  , Q_{2d}$, and $2d$ $(b, c)$-paths: $R_1^*, R_2, \ldots , R_{2d}$.

		\begin{figure}[h]
			\centering
			\includegraphics[width=0.95\columnwidth]{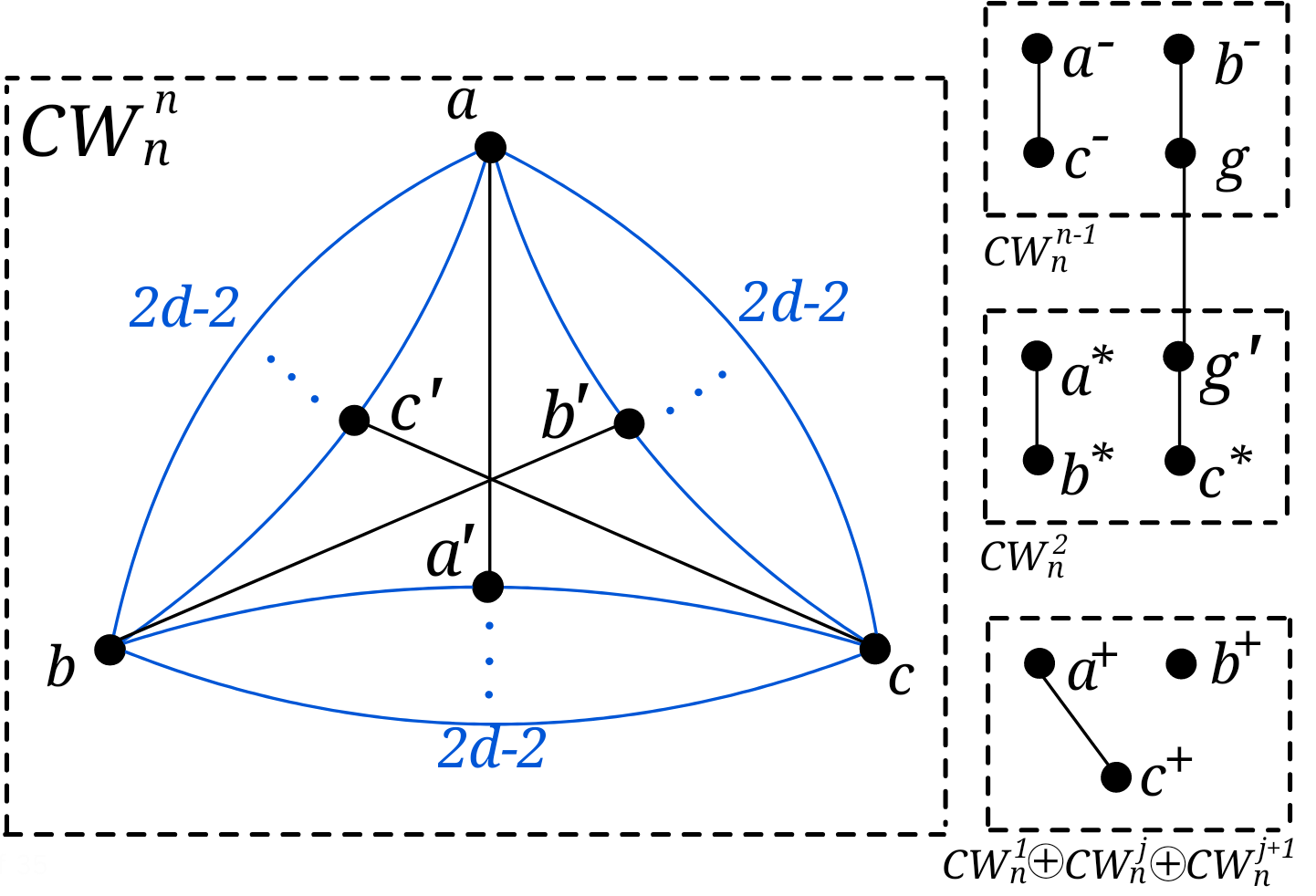} 
			\caption{Illustration of Subcase 1.2.2}
			\label{fig4}
		\end{figure}
	
		\textbf{Subcase 1.2.2.} For any $P\in L$, $d(P)=1$.

		In this case, there are three $(a, b)$-paths of length $2$, namely $aa^\prime b$, $ab^\prime b$ and $ac^\prime b$. In view of the construction of the generating set $S_2$ of $CW_n$, there must exist some $j\in[2,n{-}2]$ such that $a^{-1}b=(1~j~j{+}1)$, since  
		$(1~j~j{+}1) = (1~j)(j~j{+}1) = (1~j{+}1)(1j) = (j~j{+}1)(1~j{+}1)$.
		Similarly, there are three $(a, c)$-paths (namely $aa^\prime c, ab^\prime c, ac^\prime c$) of length $2$ and three $(b, c)$-paths (namely $ba^\prime c, bb^\prime c, bc^\prime c$) of length $2$, implying that $b^{-1}c =(1~j~j{+}1)$
		and $c^{-1}a = (1~j~j{+}1)$. Thus, $b=a(1~j~j{+}1)$ and $c=a(1~j{+}1~j).$
		W.l.o.g., assume that $a$ be the identity element of $S_n$, and thus
		$b = (1~j ~j{ + }1), c = (1~j{ + }1~j)$.

		Since each vertex  $a^\prime, b^\prime,$ and $c^\prime$ can serve as the internal vertex of  the $(a, b)$-paths, $(b, c)$-paths, and $(a, c)$-paths of length $2$, it suffices to prove that there exist four pairwise disjoint paths in $CW_n \setminus CW_n^n$ among the outside neighbors of $a, b,$ and $c$.
		
		If $j\in [3,n{-}3]$, then one can check that $a^-,b^-,c^-\in CW_n^{n-1}$, $a^*,b^*,c^*\in CW_n^2$, and $a^+,b^+,c^+\in CW_n^{1}\oplus CW_n^{j}\oplus CW_n^{j+1}$.
		By Lemma \ref{lem2.1} (iii), 
		we can find an edge $(g,g^\prime )\in E_{2,n-1}(CW_n)$ such that $g\in CW_n^{n-1}\setminus\{a^-,b^-,c^-\}$ and $ g^\prime \in CW_n^2\setminus\{a^*,b^*,c^*\}$.
		Let $X_1=\{a^-,b^-\}$ and $Y_1=\{c^-,g\}$.
		By Lemma \ref{lem2.3}, there exist two disjoint $(X_1,Y_1)$-path in $CW_n^{n-1}$ (w.l.o.g., we assume that these two paths are the $(a^-,c^-)$-path and the $(b^-,g)$-path, respectively).
		Let $X_2=\{a^*,g^\prime\}$ and $Y_2=\{b^*,c^*\}$. 
		By Lemma \ref{lem2.3}, there exist two disjoint $(X_2,Y_2)$-paths in $CW_n^2$ (w.l.o.g., we assume that these two paths are the $(a^*,b^*)$-path and the $(g^\prime,c^*)$-path, respectively).
		In addition, since $\kappa(CW_n^{1}\oplus CW_n^{j}\oplus CW_n^{j+1})\geq1$, there exist an $(a^+,c^+)$-path in $ CW_n^{1}\oplus CW_n^{j}\oplus CW_n^{j+1}$.
			Thus, referring to Figure \ref{fig4}, we have obtained four pairwise disjoint paths in $CW_n \setminus  CW_n^n$ among the outside neighbors of $a, b$, and $c$.

		If $j=2$, then $a=(1), b=(1~2~3),c=(1~3~2)$. In this case, $a^+,c^*\in CW_n^1$, $a^*,b^+\in CW_n^2$, $c^+,b^*\in CW_n^3$, and $a^-,b^-,c^-\in CW_n^{n-1}$. 
	By Lemma \ref{lem2.3}, we can easily obtain four pairwise disjoint paths in $CW_n \setminus  CW_n^n$ among the outside neighbors of $a, b$, and $c$
		 (for example: an $(a^+,c^*)$-path in $CW_n^1$, an $(a^*,b^+)$-path in $CW_n^2$, a $(c^+,b^*)$-path in $CW_n^3$, and an $(a^-,b^-)$-path in $CW_n^{n-1}$).
		
		If $j=n{-}2$, then $a=(1), b=(1~n{-}2~n{-}1),c=(1~n{-}1~n{-}2)$. In this case, $a^+,b^-\in CW_n^1$, $b^+,c^-\in CW_n^{n-2}$, $a^-,c^+\in CW_n^{n-1}$, and $a^*,b^*,c^*\in CW_n^2$. 
		Similarly, by Lemma \ref{lem2.3}, we can also obtain four pairwise disjoint paths in $CW_n \setminus  CW_n^n$ among the outside neighbors of $a, b$, and $c$.
		
		\begin{figure}[h]
			\centering
			\includegraphics[width=0.9\columnwidth]{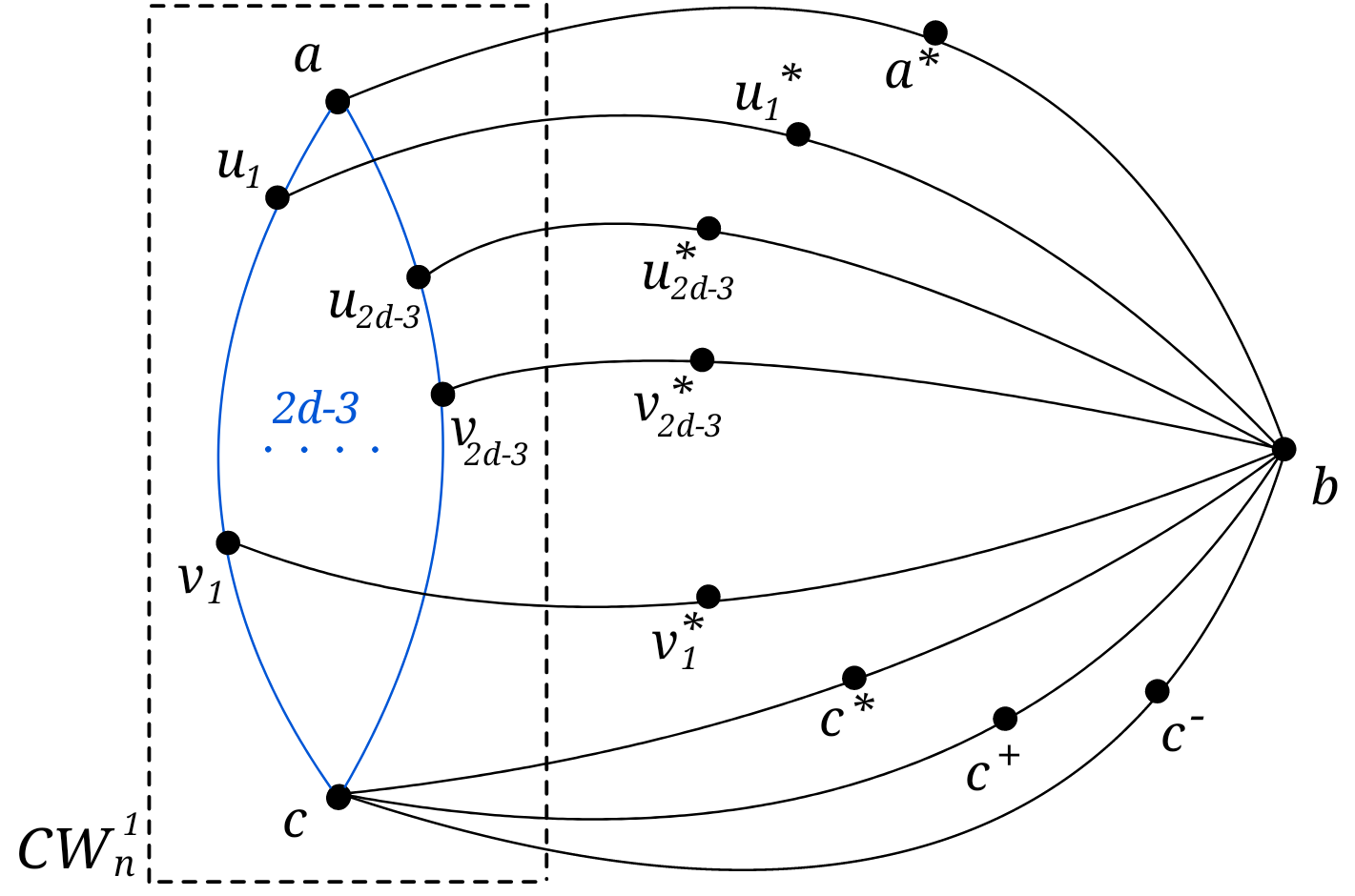} 
			\caption{Illustration of Case 2}
			\label{fig5}
		\end{figure}

		\textbf{Case 2.}  $a,b,$ and $c$ belong to two different copies.
		
		W.l.o.g., assume that $a,c\in CW_n^1$ and $b\in CW_n\setminus CW_n^1$.
		Since $CW_n^1$ is isomorphic to $BS_ {n-1}$, it follows from Lemma \ref{lem2.1} (i) that $\kappa(CW_n^1)=2n-5=4d-3$, and thus, in $CW_n^1$, there are $4d-3$ internally disjoint $(a, c)$-paths denoted by $P_i$ $(i\in [4d-3])$. By Lemma \ref{lem2.1} (ii), $|CN_{CW_n}(a,c)|\leq3$, which implies that among the $4d-3$  $(a, c)$-paths, there are at least $4d-6$ $(a, c)$-paths of length at least $3$. 
		Note that  $2d-3 \leq 4d-6$ for $d\geq 2$. W.l.o.g., assume that the length of $P_i$  for $i\in [2d-3]$  is at least $3$.
		
		For each $i\in [2d-3]$, let $N_{P_i}(a)=\{u_i\}$, $N_{P_i}(c)=\{v_i\}$, $M_a=\{u_i \mid i\in [2d-3]\}$, and $M_c=\{v_i \mid i\in [2d-3]\}$. Since $d(P_i)\geq 3$ for each $i \in [2d-3]$, 
		$M_a\cap M_c=\emptyset$. Let $N = \{u_i^* \mid i\in [2d-3]\} \cup \{v_i^* \mid i\in [2d-3]\} \cup \{a^*,c^+,c^-,c^*\}$. Then $|N|= 4d-2=2n-4$ and $N \subseteq V(CW_n\setminus CW_n^1 )$. By Lemma \ref{lem2.7}, $\kappa(CW_n\setminus CW_n^1)\geq2n-4$. Consequently, referring to Figure \ref{fig5}, Lemma \ref{lem2.2} implies the existence of $4d-2$ internally disjoint $(b, N)$-paths in $CW_n\setminus CW_n^1$,  which in turn yield $2d-2$ $(a, b)$-paths and $2d$ $(b, c)$-paths. Combining these paths with $2d$ $(a,c)$-path $P_i$ $(i\in [2d-2,4d-3])$,  the structure illustrated in Figure \ref{fig2} (b) is obtained.

		\textbf{Case 3.}
		$a,b,$ and $c$ belong to three different copies, respectively.
		
		W.l.o.g., assume that $a\in CW_n^1, b\in CW_n^2$, and $c\in CW_n^3$.
		
		Let $W_1=\{w \,|\, w\in V(CW_n^1)\setminus\{a,b^+,b^-,b^*,c^+,c^-,c^*\},w^*\in  V(CW_n^2)\}$,
		
		~~~~~$W_2=\{w \,|\, w\in V(CW_n^1)\setminus\{a,b^+,b^-,b^*,c^+,c^-,c^*\},w^*\in  V(CW_n^3)\}$,
		
		~~~~~$W_3=\{w \,|\, w\in V(CW_n^2)\setminus\{b,a^+,a^-,a^*,c^+,c^-,c^*\},w^*\in  V(CW_n^3)\}$.
		
		Clearly, the sets $W_1, W_2,$ and $W_3$ are pairwise disjoint, with $|W_i|\geq (n-2)!-3\geq 2d-1$ for each $i=1,2,3$.
		
		Let $\widehat{CW_n}=CW_n\setminus(CW_n^1\oplus CW_n^2\oplus CW_n^3)$. 
		By Lemma \ref{lem2.1} (iv), each of the vertices $a,b$, and $c$ has at least one and at most three outside neighbor in $\widehat{CW_n}$.
	 W.l.o.g., assume that $a^*,b^*,c^*\in \widehat{CW_n}$. 
		Next, we shall complete the proof of Case 3 by considering the following three exhaustive subcases:
		\begin{figure}[h]
			\centering
			\includegraphics[width=0.67\columnwidth]{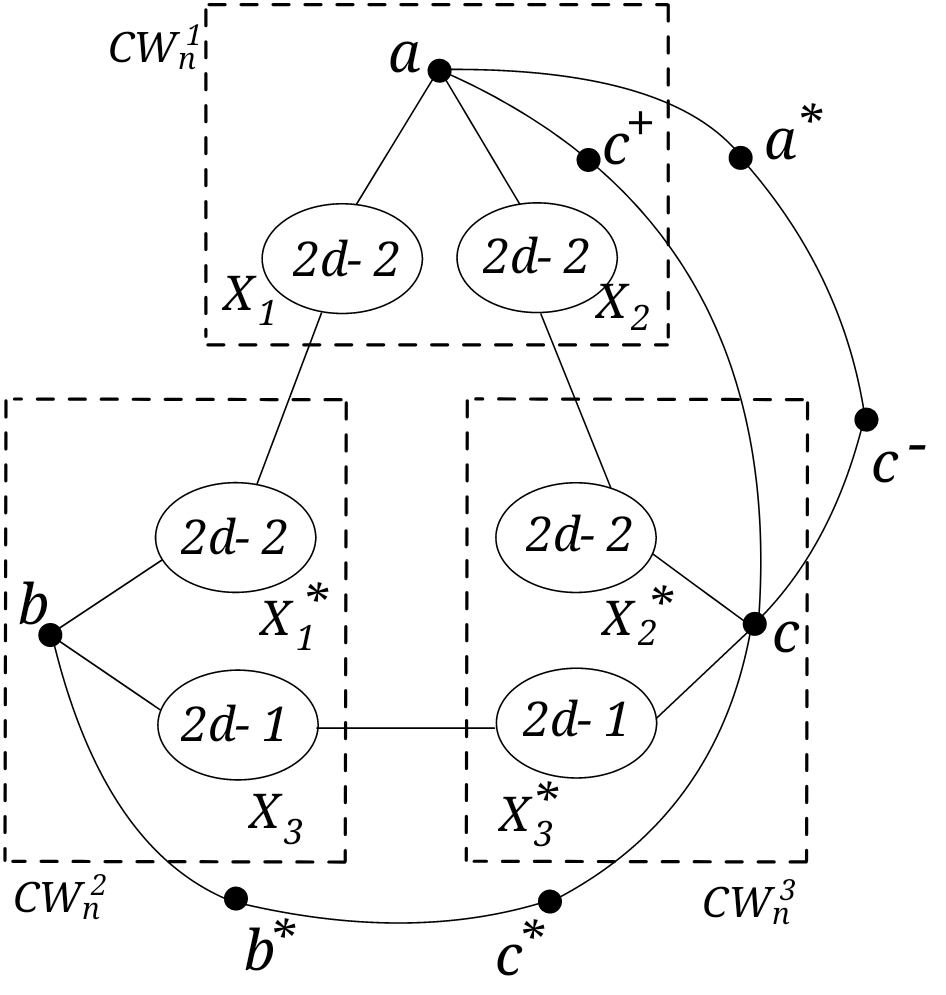} 
			\caption{Illustration of Subcase3.1}
			\label{fig6}
		\end{figure}
	
		\textbf{Subcase 3.1.}  
		There exists at least one vertex in $\Omega$ that has exactly two outside neighbors in $\widehat{CW_n}$.
		
	W.l.o.g., assume that $c^*,c^-\in \widehat{CW_n}$ and $c^+\in {CW_n^1}$.
	Since $|W_i|\geq 2d-1$ for each $i =1,2,3$, 
	we can find $X_i \subseteq W_i$ with $|X_i| = 2d-2$ for $i = 1, 2$ and $X_3 \subseteq W_3$ with $|X_3| = 2d - 1$.
	Let  $X_i^* = \{x^* \mid x \in X_i\}$ for $i = 1,2,3$. 
		Obviously, $(X_1 \cup X_2\cup\{c^+\}) \subseteq V(CW_n^1)$, $(X_1^* \cup X_3) \subseteq V(CW_n^2)$, and $(X_2^* \cup X_3^*) \subseteq V(CW_n^3)$. 
		Since  $|X_1 \cup X_2\cup\{c^+\}| = |X_1^* \cup X_3|=|X_2^* \cup X_3^*| =4d-3=\kappa(CW_n^i)$ for $i=1,2,3$, it follows from Lemma \ref{lem2.2} that there exist $4d-3$ internally disjoint $(a, X_1 \cup X_2\cup\{c^+\})$-paths in $CW_n^1$, $4d-3$ internally disjoint $(b, X_1^* \cup X_3)$-paths in $CW_n^2$, and $4d-3$ internally disjoint $(c, X_2^* \cup X_3^*)$-paths in $CW_n^3$. Note that there is a perfect matching between $X_i$ and $X_i^*$ for  $i\in [3]$. Referring to Figure \ref{fig6}, we can obtain $2d-2$ internally disjoint $(a,b)$-paths: $a-X_1-
		X_1^*-b$, $2d-1$ internally disjoint $(a,c)$-paths: $a-X_2-
		X_2^*-c$ and $a-c^+c$, and $2d-1$ internally disjoint $(b,c)$-paths: $b-X_3-X_3^*-c$.
		
			Let $X=\{c^-,c^*\}$ and $Y=\{a^*,b^*\}$. Then $X, Y\subseteq V(\widehat{CWn})$.
			By Lemma \ref{lem2.6}, $\kappa(\widehat{CW_n})\geq2n-5\geq2$, and thus by Lemma \ref{lem2.3}, there exist two disjoint $(X,Y)$-paths in $\widehat{CW_n}$ (w.l.o.g., we assume that these two paths are the $(a^*,c^-)$-path and the $(b^*,c^*)$-path, respectively). Thus, referring to Figure \ref{fig6}, the structure illustrated in Figure \ref{fig2} (b) is obtained.

			\begin{figure}[h]
			\centering
			\includegraphics[width=0.7\columnwidth]{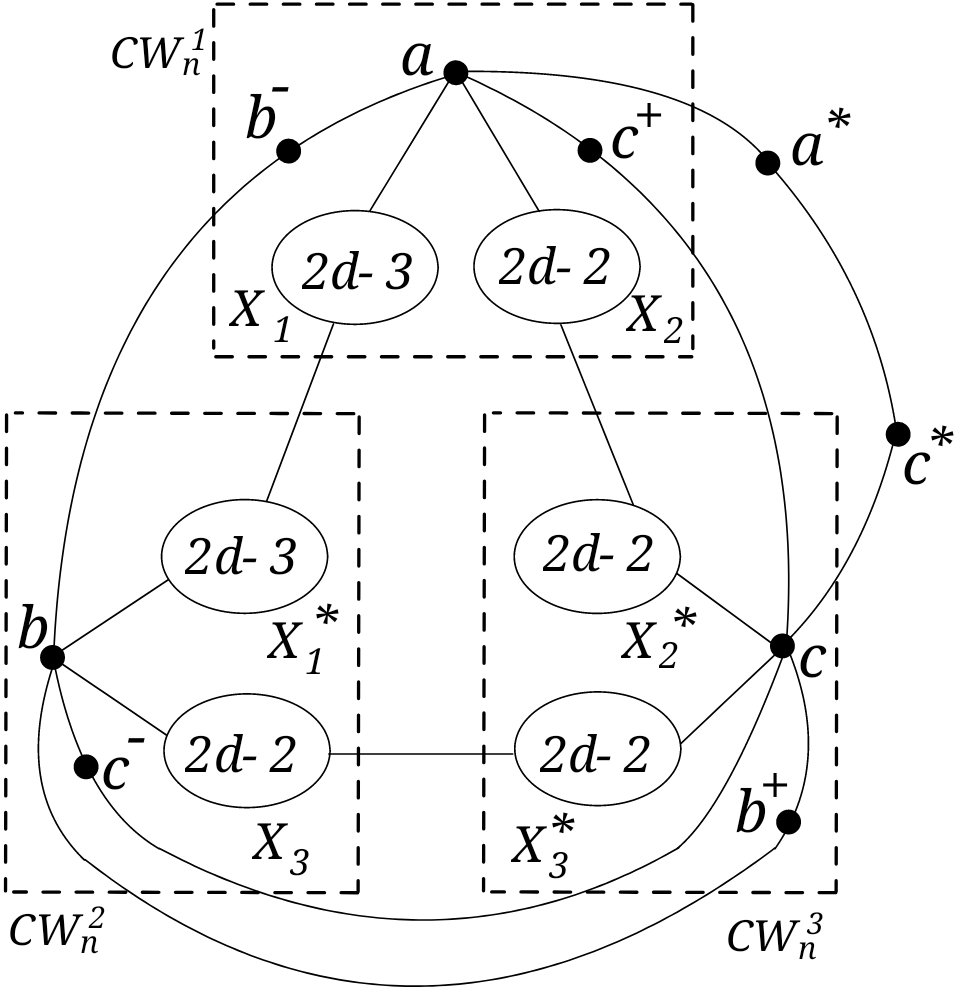} 
			\caption{Illustration of Subcase3.2}
			\label{fig7}
		\end{figure}
	
		\textbf{Subcase 3.2.} 
		There exist at least two vertices in $\Omega$ each of which has exactly one outside neighbors in $\widehat{CW_n}$.
		
		W.l.o.g., assume that $b^*,c^*\in \widehat{CW_n}$, $c^+,b^-\in {CW_n^1}$, $c^-\in {CW_n^2}$, and $b^+\in {CW_n^3}$. In much the
		same way as Subcase 3.1,
		 we can find $X_1 \subseteq W_1$ with $|X_1| = 2d - 3$ and $X_i \subseteq W_i$ with $|X_i| = 2d-2$ for $i = 2,3$. 
		 Let $X_i^* = \{x^* \mid x \in X_i\}$. Similarly, by Lemma \ref{lem2.2}, there exist $4d-3$ internally disjoint $(a,X_1\cup X_2\cup\{b^-,c^+\})$-paths in $CW_n^1$, $4d-4$ internally disjoint $(b,X_1^*\cup X_3\cup\{c^-\})$-paths in $CW_n^2$, and $4d-3$ internally disjoint $(c,X_2^*\cup X_3^*\cup\{b^+\})$-paths in $CW_n^3$. In addition, since $\kappa(\widehat {CW_n})\geq 1$, there is an $(a^*,c^*)$-path in $\widehat {CW_n}$. Thus, referring to Figure \ref{fig7}, we can obtain $2d-2$ internally disjoint $(a,b)$-paths (namely $a-X_1-X_1^*-b$ and $a-b^{-}b$),  $2d$ internally disjoint $(a,c)$-paths (namely $a-X_2-X_2^*-c$, $aa^*-c^*c$, and $a-c^+c$), and $2d$ internally disjoint $(b,c)$-paths (namely $b-X_3-X_3^*-c$, $b-c^{-}c$, and $bb^+-c$).
	 
			\begin{figure}[h]
			\centering
			\includegraphics[width=0.7\columnwidth]{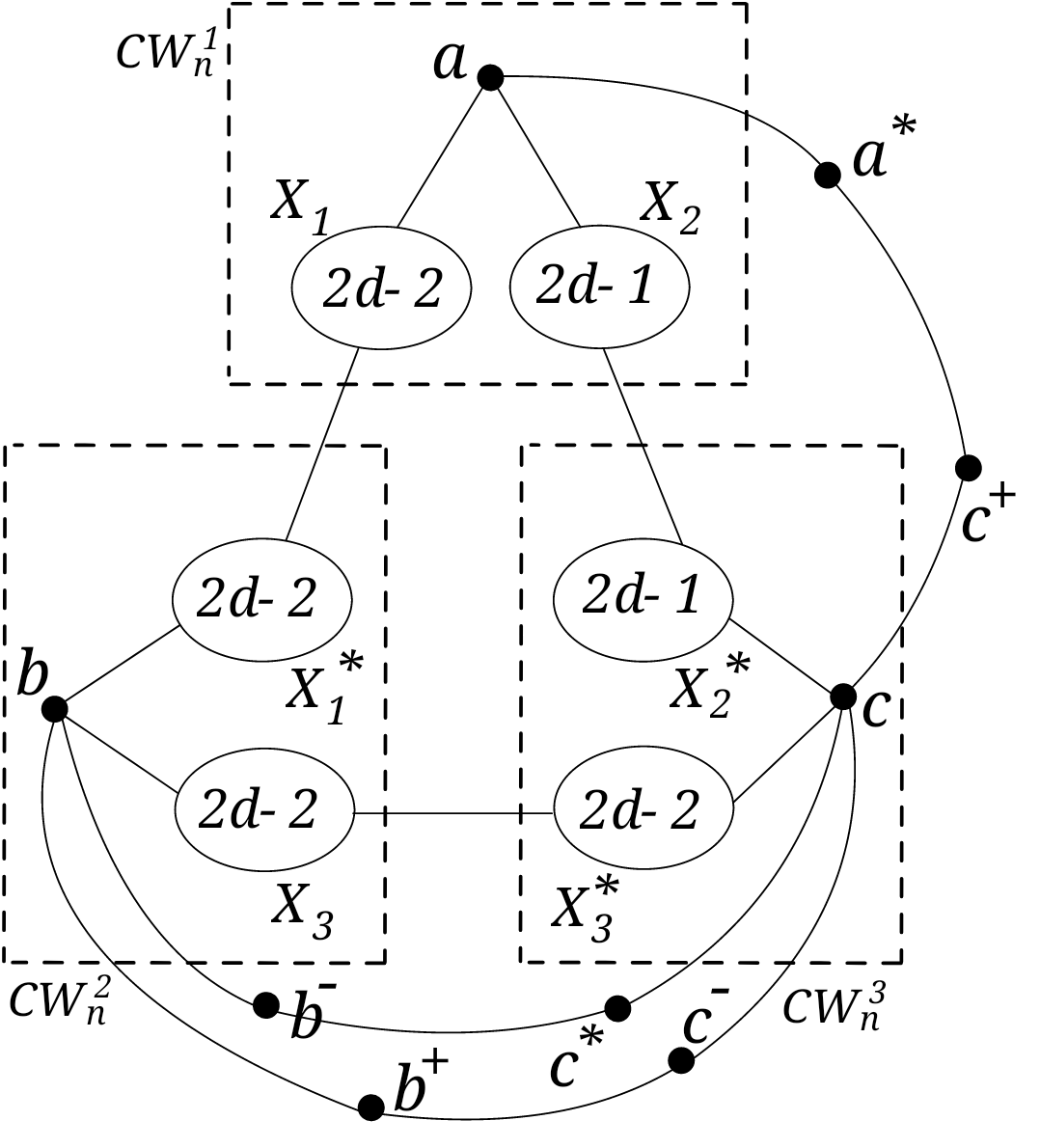} 
			\caption{Illustration of Subcase3.3}
			\label{fig8}
		\end{figure}
		
		\textbf{Subcase 3.3.} 
		There exist at least two vertices in $\Omega$ each of which has exactly three outside neighbors in $\widehat{CW_n}$.
		
		W.l.o.g., assume that  $a^*,b^+,b^-,b^*,c^+,c^-,c^*\in \widehat{CW_n}$. In much the
		same way as Subcase 3.1, 
		we can find $X_2 \subseteq W_2$ with $|X_2| = 2d - 1$ and $X_i \subseteq W_i$ with $|X_i| = 2d-2$ for $i = 1,3$. Let $X_i^* = \{x^* \mid x \in X_i\}$. Similarly, by Lemma \ref{lem2.2}, there exist $4d-3$ internally disjoint $(a,X_1\cup X_2)$-paths in $CW_n^1$, $4d-4$ internally disjoint $(b,X_1^*\cup X_3)$-paths in $CW_n^2$, and $4d-3$ internally disjoint $(c,X_2^*\cup X_3^*)$-paths in $CW_n^3$.
		 Let $X=\{c^+,c^-,c^*\}$ and $Y=\{a^*,b^+,b^-\}$. Then $X,Y \subseteq V(\widehat{CW_n})$. Since $\kappa(\widehat {CW_n})\geq 3$, it follows from Lemma \ref{lem2.3} that there exist three disjoint $(X,Y)$-paths in $\widehat{CW_n}$ (w.l.o.g., we assume that these three paths are the $(a^*,c^+)$-path, the $(b^-,c^*)$-path, and the $(b^+,c^-)$-path, respectively). Thus, referring to Figure \ref{fig8}, we can obtain $2d-2$ internally disjoint $(a,b)$-paths (namely $a-X_1-X_1^*-b)$,  $2d$ internally disjoint $(a,c)$-paths(namely $a-X_2-X_2^*-c$, $aa^*-c^{+}c$), and $2d$ internally disjoint $(b,c)$-paths (namely $b-X_3-X_3^*-c$, $bb^{-}-c^*c$, and $bb^+-c^{-}c)$. $\square$
		
	\end{theorem}

	\section{Main result}

	\begin{theorem}
		For $n\geq4$, $\pi_3(CW_n)=\lfloor\frac{6n-9}4\rfloor$.
		\begin{figure}[htbp]
		\centering
		\subfloat[$n=2d\,(d\geq2)$]{
			\includegraphics[width=0.48\columnwidth]{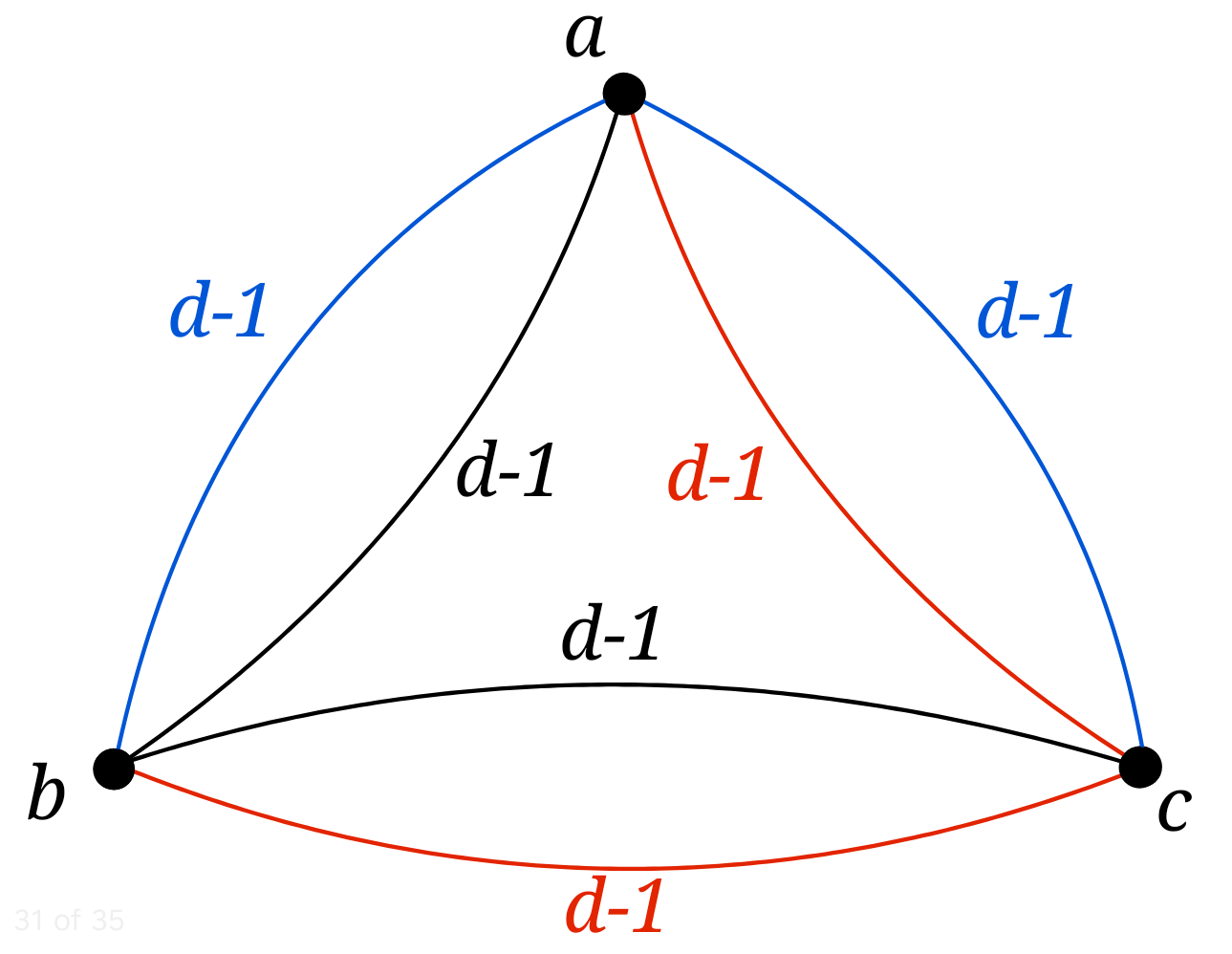}
			
		}
		\hfill
		\subfloat[$n=2d+1\,(d\geq2)$]{
			\includegraphics[width=0.48\columnwidth]{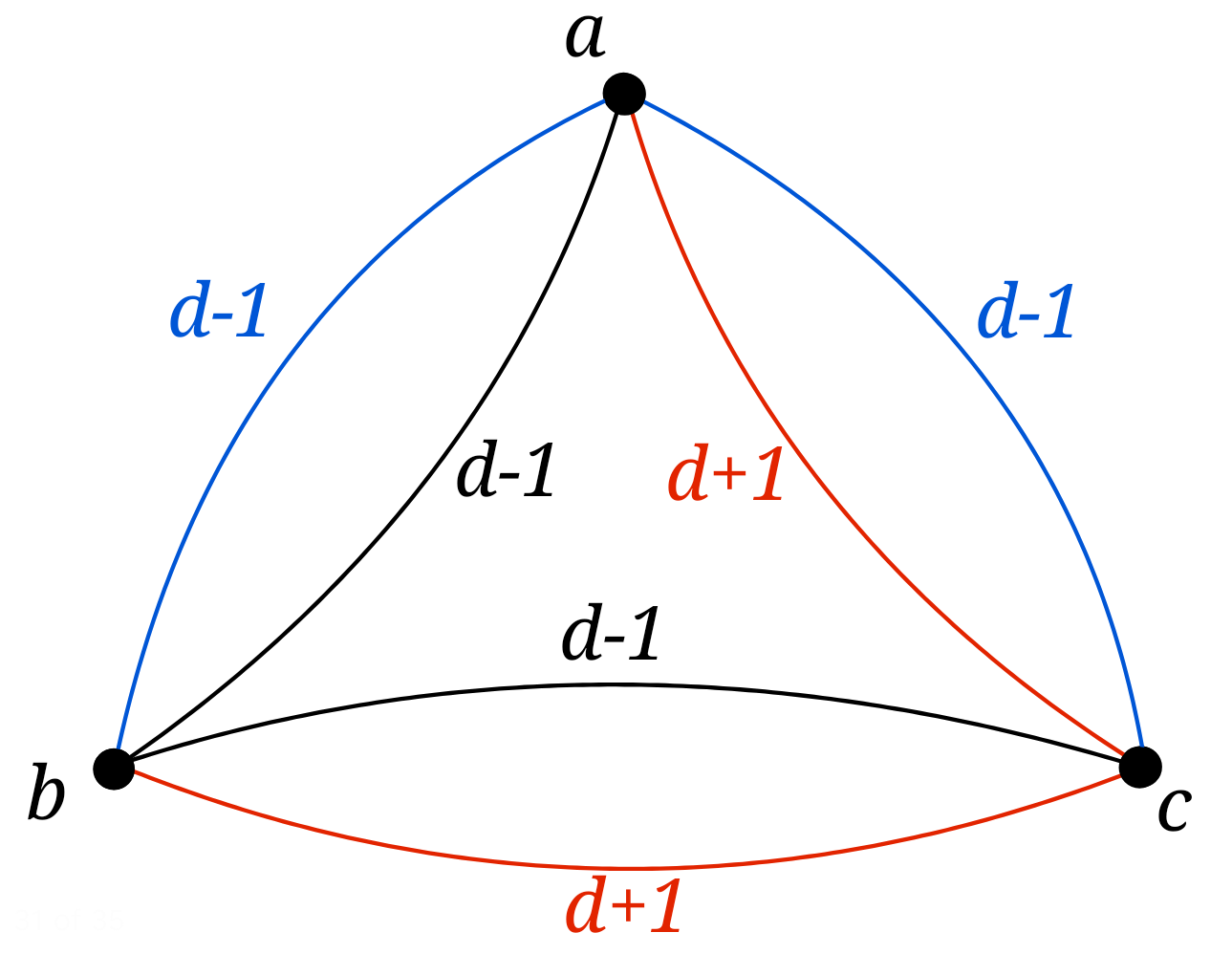}
			
		}
		\caption{Illustrations of Theorem 4.1}
		\label{fig9}
	\end{figure}

		\noindent\textbf{Proof.}
		Since $CW_n$ is a $(2n-2)$-regular connected graph, it follows from Lemmas \ref{lem2.5} and \ref{lem2.9} that   $\pi_3(CW_n)\leq\lfloor\frac{3(2n-2)-3}4\rfloor=\lfloor\frac{6n-9}4\rfloor$.
		Let $\Omega=\{a,b,c\}$ be any $3$-subset of $V(CW_n)$.
		Next, we show that $\pi_3(CW_n)\geq\lfloor\frac{6n-9}4\rfloor$ by considering the parity of $n$.

		\textbf{Case 1.} $n=2d\,(d\geq2)$.
		
		By Theorem \ref{the3.1},
		$CW_n$ contains the structure illustrated in Figure \ref{fig2} (a). 
		Referring to Figure \ref{fig9} (a), 
		we first pair $d{-}1~(a,b)$-paths with $d{-}1~(a,c)$-paths. 
		Then, pair the remaining $d{-}1~(a,b)$-paths with $d{-}1~(b,c)$-paths. 
		Finally, pair the remaining $d{-}1~(b,c)$-paths with the remaining $d{-}1~(a,c)$-paths.	
		In this way, we can find $3d-3$ internally disjoint $\Omega$-paths, and thus $\pi_3(CW_n)\geq3d-3=\lfloor\frac{6n-9}4\rfloor$.
		
		\textbf{Case 2.} $n=2d+1\,(d\geq2)$.
		
		By Theorem \ref{the3.1},
		$CW_n$ contains the structure illustrated in Figure \ref{fig2} (b). 
		Referring to Figure \ref{fig9} (b), 
		we first pair $d{-}1~(a,b)$-paths with $d{-}1~(a,c)$-paths. 
		Then, pair the remaining $d{-}1~(a,b)$-paths with $d{-}1~(b,c)$-paths. 
		Finally, pair the remaining $d{+}1~(b,c)$-paths with the remaining $d{+}1~(a,c)$-paths.	
		In this way, we can find $3d-1$ internally disjoint $\Omega$-paths, and thus $\pi_3(CW_n)\geq3d-1=\lfloor\frac{6n-9}4\rfloor$.
		
		Thus, we have proved that $\pi_3(CW_n)=\lfloor\frac{6n-9}4\rfloor$ for $n\geq4$. $\square$

	\end{theorem}

	\end{document}